\documentclass[twoside,reqno,10pt]{amsart}

	 \newcommand{\fracvar}[4]{\ensuremath{ \upsilon_{ #3 }^{ #4} [ #1 ] \left(  #2 \right)   }}
	 \newcommand{\fracvarplus}[3]{ \fracvar {#1}{#2}{#3}{\epsilon+} }
	 \newcommand{\fracvarmin}[3]{ \fracvar {#1}{#2}{#3}{\epsilon -} }
	\newcommand{\llim}[3]{\ensuremath{ \lim\limits_{ #1 \rightarrow #2} #3 }}
	\newcommand{\fclass}[2]{\ensuremath{  \mathbb{#1}^{\, #2} }}
	
		\newcommand{\holder}[1]{\fclass{H}{#1} }

 \newtheorem{theorem}{Theorem}
 \newtheorem{lemma}{Lemma}
 \newtheorem{corollary}{Corollary}
 \newtheorem{proposition}{Proposition}
\newtheorem{definition}{Definition}
 \newtheorem{remark}{Remark}
 \newtheorem{example}{Example}

\usepackage{graphicx,color,psfrag}

 \usepackage{lineno}

\title[Fractal variation of functions]
{Fractional variation of H\"olderian functions}
  
\author {Dimiter Prodanov}
\address{Correspondence: Environment, Health and Safety, IMEC vzw, Kapeldreef 75, 3001 Leuven, Belgium}

\begin{document}

%
%
%
%

\maketitle
\section*{Preamble}
\large{Please cite this paper as \\
D. Prodanov, "Fractional variation of H\"olderian functions", 
Fract. Calc. Appl. Anal., Vol. 18, No 3 (2015), pp. 580-602 \\
DOI: 10.1515/fca-2015-0036}

Proofs of the results are given in the full paper or available upon request.
This is a shortened version of the paper in compliance with the copyright statement. 

\section{Introduction}
\label{seq:intro}
 
Fractional derivatives and fractional calculus have long history since the time of H\^opital and Leibniz \cite{Oldham1970,Ross1977}. 
However, only relatively recently fractional calculus has been recognized as a tool for modeling physical and biological problems (see for example \cite{Caputo1967,Caputo1971,Bagley1983}). 
Some classical definitions of fractional derivatives (for example by Riemann and Liouville) are based on extension of the Cauchy integral into non-integer order. 
However, such derivatives are difficult to compute and their geometric interpretation is unclear because of their non-local character.
In particular, there is no relationship between the local geometry of the graph of function and its fractional derivative \cite{Adda1997, Adda2005}. A definite disadvantage of the Riemann-Liouville approach is that the fractional derivative of a constant is not zero. 
This was the starting point for the modified definitions of Caputo \cite{Caputo1967} and Jumarie \cite{Jumarie2006}.

Recently, natural science applications have also inspired the development of local fractional derivatives \cite{Kolwankar1997a, Adda2001}. 
The theory of such derivatives is still in its infancy and there are few available results   \cite{Kolwankar1997a, Adda2001,  Babakhani2002, Chen2010}. 
The starting point of Kolwankar and Gangal \cite{Kolwankar1997a} was the Riemann-Liouville approach (recent review in \cite{Kolwankar2013}).
 On the other hand, Ben Adda and Cresson  \cite{Adda2001} introduced from the start a difference operator based definition easily transferable to integer-ordered derivatives. 
The correspondences between the integral approach and the quotient difference approach have been further investigated by  Chen et al. 
\cite{Chen2010} and some of the initial results have been clarified and corrected in \cite{Adda2013}.

In this paper I present a method based on the so-called \textit{fractal variation operators}, which can provide a simple way of local characterization of singular and scaling behavior of continuous functions.
Fractal variation operators are constructed from power scaling of finite difference operators. 
Application of this approach can be especially suitable for characterization of  H\"olderian (especially non-differentiable) functions.
One of the main results comprises the calculation of the fractal variation of Cauchy sequences leading to the Dirac's $\delta$-function. 

The manuscript is organized as follows. 
Section \ref{sec:definitions} introduces the general definitions and notations. 
Section \ref{sec:holderprop} introduces H\"olderian functions and demonstrates their properties used in further proofs.
Section \ref{sec:frdiff} introduces the \textit{fractal variation} operators.
Sections \ref{sec:apps} demonstrates some of applications of \textit{fractal variation} to smooth and singular functions.
The fractal variation of a function in the infinitesimal limit corresponds to the definition of local fractional derivative introduced by Ben Adda and Cresson \cite{Adda2001, Adda2005, Chen2010}.
Notably, these authors define the local factional derivative as 
\[
 \frac{d^\alpha}{d x^\alpha}f(x) := \Gamma (1+\alpha) \llim{\epsilon}{0}{\frac{ f(x+ \epsilon) - f(x) }{\epsilon ^\alpha}}
\]
However, according to the main results in the present work (Theorems  and \ref{th:singul}) this derivative has very few non-infinitesimal or non-divergent values. 
Therefore, I prefer the term "variation" over a  derivative.
Moreover, such operators can be more useful in finite difference settings, e.g. for numerical applications.

\section{General definitions and notations}
\label{sec:definitions}

Along the text I consistently use square brackets for the arguments of operators and round brackets for the arguments of functions. 
The term \textit{function}  denotes the mapping $ f: \mathbb{R} \mapsto \mathbb{R} $. The notation $f(x)$ is used to refer to the value of the function at the point \textit{x}. 
By \fclass{C}{0} is denoted the class of functions that are continuous and by \fclass{C}{n} the class of \textit{n}-times differentiable functions where  $n  \in \mathbb{N}$.
$Dom[ f ]$ denotes the domain of definition of the function $f(x)$.

\section{H\"olderian functions}
\label{sec:holderprop}
\begin{definition}
\label{def:holder}
Let \holder{\alpha} be the class of H\"older   $\mathbb{C}^0$  functions of degree $\alpha$, $\alpha \in (0,\, 1)$,
That is, $\forall f (x) \in \holder{\alpha} $ there exist two positive constants 
$C, \delta \in \mathbb{R} $ for $  x,y \in Dom[ f ]$ such that for $|x-y| \leq \delta$ the following inequality holds
\[
| f (x) - f (y) |  \leq C |x-y|^\alpha 
\]

Following Mallat and Hwang \cite{Mallat1992} the definition can be extended to orders greater than one in the following way. 
Let \holder{n+ \alpha} be the class of \fclass{\, C}{0}  double H\"older functions of degree $n+\alpha$ for which
\[
| f (x) - f (y) - P_n (x-y) |  \leq C |x-y|^{n +\alpha} 
\]
where $P_n  $ is a real-valued polynomial of degree $n \in \fclass{N}{}$ of the form
\[
  P_n (z):= \sum\limits_{k=1}^{n}{ a_k z^k} 
\]
where $P_0(z)=0$ and $\alpha \in (0,\, 1)$.
\end{definition}

Under this definition we will focus mainly on functions for which $0< \alpha<1$. These functions will be further called by the term \textbf{H\"olderian}.
\begin{definition}
	\label{def:deltas}
Let the difference parametrized operators acting on a function $f(x)$ be defined in the following way
\begin{align}
  \Delta^{+}_{\epsilon} [f](x) & :=  f(x + \epsilon) - f(x) \\
   \Delta^{-}_{\epsilon} [f](x) & :=  f(x) - f(x - \epsilon)  \\
    \Delta^{2}_{\epsilon} [f](x) &:=  f(x + \epsilon) -2 f(x) + f(x - \epsilon)
\end{align}
where $\epsilon>0$. The first one we refer to as \textit{forward difference} operator, 
the second one we refer to as \textit{backward difference} operator and the third one as
\textit{2\textsuperscript{nd} order difference} operator.
\end{definition}
\begin{lemma}[Difference composition lemma]
\label{th:diffcomp}
The 2\textsuperscript{nd} order difference operator is a composition of the backward and forward difference operators.
\[
  \Delta^{2}_{\epsilon} =  \Delta^{+}_{\epsilon} \circ \Delta^{-}_{\epsilon} = \Delta^{-}_{\epsilon} \circ \Delta^{+}_{\epsilon} = \Delta^{+}_{\epsilon} - \Delta^{-}_{\epsilon} 
  \]
\end{lemma}


\begin{theorem}
\label{prop:holdern}
Let $f(x) \in \holder{n+\alpha}$ in the interval $[x \ x+ \epsilon]$, where $n \in \fclass{N}{}$ is a natural number and $ \alpha \in (0 \ 1)$. 
Then $ f(x) \in \fclass{C}{n}$ in the interval $[x \ x+ \epsilon]$.
\end{theorem}
\begin{corollary}
If 
$
 | f (x) - f (y) - P_n (x-y) |  \leq C^\prime |x-y|^{n +\alpha} 
$
for $n  \in \fclass{N}{}$ and $0<\alpha <1$
then
\[
  P_n (x-y)= \sum\limits_{k=1}^{n}{ \frac{1}{k!}  f^{(k)} (x) \,(x-y)^k} 
\]
\end{corollary}
\section{Fractal variation operators}
\label{sec:frdiff}

\begin{definition}
\label{def:fracvar}
Let the \textit{Fractal Variation} operators be defined as
\begin{align}
\label{eq:fracvar1}
	 \fracvarplus {f}{x}{\beta} := \frac{\Delta^{+}_{\epsilon} [f] (x) }{\epsilon ^\beta}  =\frac{ f(x+ \epsilon) - f(x) }{\epsilon ^\beta} 
	  \\
	 \fracvarmin {f}{x}{\beta} :=  \frac{\Delta^{-}_{\epsilon} [f] (x) }{\epsilon ^\beta}  =\frac{ f(x)- f( x- \epsilon)  }{\epsilon  ^\beta} 
\end{align}
where  $\epsilon >0$ and $0 < \beta \leq 1 $ are real parameters and $f(x)$ is a function.
\end{definition}

It is easy to check that the Fractal variation operators  are R-linear. 
This can be formulated in the following Lemma
\begin{lemma}
\label{th:linearity}
Let \textit{K} and \textit{M} be real constants. Then
\[
 \upsilon_{\beta}^{ \epsilon +} [ K \, f (x) + M \, g(x)]  =  K \, \upsilon_{\beta}^{ \epsilon +} [ f] (x) +  M \, \upsilon_{\beta}^{ \epsilon +} [ f] (x)
\] 
and
\[
 \upsilon_{\beta}^{ \epsilon -} [ K \, f (x) + M \, g(x)]  =  K \, \upsilon_{\beta}^{ \epsilon -} [ f] (x) +  M \, \upsilon_{\beta}^{ \epsilon -} [ f] (x)
\]
\end{lemma}

\begin{definition}
	\label{def:transl}
	Let the translation operator acting on the function $f(x)$ be defined as
	\[
	\mathcal{T}_{\pm \epsilon} [f] (x) :=f( x \pm \epsilon)
	\]
	where $ \epsilon$ is a real positive parameter.
\end{definition}
It is easy to check that for continuous functions the reflection of a translation about the origin is its inverse operation:
$
\mathcal{T}_{-\epsilon}= \mathcal{T}^{-1}_{\epsilon}
$.
Fractional variation is translationally invariant. This can be formulated in the following result.

\begin{proposition}
Fractal variation commutes with the translation operator in the neighborhood of a point if the argument function is defined in the neighborhood.

Let $f(x) \in \holder{\alpha}$. 
Then 
\[
\{ \upsilon_{\beta}^{ \epsilon \pm} \mathcal{T}_\epsilon - \mathcal{T}_\epsilon \upsilon_{\beta}^{ \epsilon \pm} \} [f]  = 0 [f]
\]
\end{proposition}
The proof follows directly from the  commutativity of addition in $\mathbb{R}$.
\begin{remark}
	The following relation can be established from the definitions:
	\begin{align}
	\fracvarplus {f}{x}{\beta} 
	=  \epsilon^{-\beta} \left( \mathcal{T}_{\epsilon} - \mathcal{I} \right) [f] (x) \\
	\fracvarmin {f}{x}{\beta} 
	=  -\epsilon^{-\beta} \left( \mathcal{T}_{-\epsilon} - \mathcal{I} \right) [f] (x)
	\end{align}
	where $\mathcal{I}$ is the identity operator.
\end{remark}
\begin{proposition}
	\label{prop:map}
The right and left fractal variation operators are mapped to each other by translation assuming that the argument function is defined in the domain of the operators. 
\end{proposition}
\begin{theorem}[Duality of the limit of variation]
	\label{th:dualvar}
	If $f(x)$ is defined in $[ x - \epsilon, x+\epsilon]$ and the limit 
	\llim{\epsilon}{0}{\fracvarplus{f}{x}{\beta} } exists and 
	$
	\llim{\epsilon}{0}{\dfrac{\Delta_{\epsilon} ^2 [f] (x)}{\epsilon^\beta} } = 0 
	$
	then 	$
	\llim{\epsilon}{0}{\fracvarplus{f}{x}{\beta} } =  \llim{\epsilon}{0}{\fracvarmin{f}{x}{\beta} } 
	$.
\end{theorem}
\begin{corollary}
	If \llim{\epsilon}{0}{\fracvarplus{f}{x}{\beta} } is continuous about $x$ then $\llim{\epsilon}{0}{\fracvarplus{f}{x}{\beta} } =  \llim{\epsilon}{0}{\fracvarmin{f}{x}{\beta} } $.
	If \llim{\epsilon}{0}{\fracvarmin{f}{x}{\beta} } is continuous about $x$ then $\llim{\epsilon}{0}{\fracvarplus{f}{x}{\beta} } =  \llim{\epsilon}{0}{\fracvarmin{f}{x}{\beta} } $.
\end{corollary}
\begin{lemma}
	\label{th:varmin}
Let $f(x) \in \holder{\alpha}$. 
Then $ \fracvarplus{f}{x}{\beta}  \in \holder{\alpha -\beta}$ and $ \fracvarmin{f}{x}{\beta}  \in \holder{\alpha -\beta}$.
\end{lemma}
\begin{theorem}[Compound variation theorem]
\label{th:compvar1}

The fractal variation of a compound function $ f\left(y(x) \right)$ can be expressed as
\begin{equation}
\label{eq:compvar3}
\fracvarplus{f (y)}{x}{\beta}  =\upsilon_{\alpha}^{\eta +}[f](y) \, \left( \fracvarplus{y}{x}{1} \right)^\alpha \epsilon^{\alpha-\beta}
\end{equation}
where $\eta= y(x+\epsilon) - y(x) \neq 0$, $\epsilon>0$ and  $\alpha>0$.
The argument function is interpreted as a variable. 
\end{theorem}
So-stated theorem can be specialized in two corollaries that are important for applications.
\begin{corollary}[First differential form]
For a compound function $ f\left(y(x) \right)$ 
\begin{equation}
\label{eq:compvar1}
\fracvarplus{f(y)}{x}{\beta} = \upsilon_{1}^{\eta +}[f](y)  \ \fracvarplus{y}{x}{\beta}
\end{equation}
where $\eta= y(x+\epsilon) - y(x)$.
\end{corollary}
\begin{corollary} [Second differential form]
For a compound function $ f\left(y(x) \right)$:
\begin{equation}
\label{eq:compvar2}
\fracvarplus{f(y)}{x}{\beta} =\upsilon_{\beta}^{\eta +}[f](y) \, \left( \fracvarplus{y}{x}{1} \right)^\beta
\end{equation}
where $\eta= y(x+\epsilon) - y(x)$.
\end{corollary}
\begin{remark}
These corollaries lead to the following limiting behavior for $\alpha<1$ for a monotonous substitution function $y(x)$
\[
\llim{\epsilon}{0} {\fracvarplus{f(y)}{x}{\beta}} =
f^\prime (y) \llim{\epsilon}{0} {\fracvarplus{y}{x}{\beta}}
\]
provided that $f^\prime (y)$ exists at $y$
and
\[
\llim{\epsilon}{0} {\fracvarplus{f(y)}{x}{\beta}} =
\left( y^\prime (x)  \right)^\beta \, \llim{\eta}{0} {\upsilon_{\beta}^{\eta +}[f](y)}
\]
provided that $y^\prime (x)$ exists at $x$.
In both formulas the argument function is interpreted as a variable. 
\end{remark}

\section{Applications}
\label{sec:apps}
\subsection{Fractal variation of the H\"older exponent}
\label{seq:limholder}

\begin{definition}
Let the sign operator acting on a function $f(x)$ be defined as
\[
sign[f](x):=\left\{
\begin{array}{ll}
	+1 ,&  f(x) \geq 0 \\
	-1 ,&  f(x)<0
\end{array} \right.
\]
\end{definition}

In the next section  the limit behavior of so-defined operators for the limiting H\"older exponent function will be demonstrated.

\begin{theorem}
\label{th:fvar}
Let $g(x)= |x|^\alpha, \ \alpha >0$. 
Then the limiting behavior of the Fractal variation $\llim {\epsilon} {0} {\upsilon_{\beta}^{ \epsilon  } [g] (x)} $ can  be summarized in   Table \ref{tab:fracvar}:

\begin{table}[!ht]
	\centering
		\begin{tabular}[t]{ll|r|r}
		\hline
		\multicolumn{2}{c}{-}  \vline   	& $x=0$   	& 	 	$ |x| >0 $  \\
		\hline
$ \alpha = \beta$ 	& $\alpha < 1$	& 1						&			 0     	 		\\
$ \alpha = \beta$ 	& $\alpha > 1$	& 1						&			 $\infty$   \\
$ \alpha > \beta$ 	& $\beta > 1$	& 0						&			 $\infty$   \\
$ \alpha < \beta$ 	& $\beta < 1$	& $\infty$				&			 $0$   \\
	 \cline{3 - 4}
	  $ \alpha <1 $	& $\beta >  1$	& 		\multicolumn{2}{c}{$\infty$}  	 		\\ 	
\multicolumn{2}{l}{$ \alpha = \beta =1 $} \vline	&  \multicolumn{2}{c}{1} \\
$ \alpha > \beta$	& $\beta < 1$	&  \multicolumn{2}{c}{ 0}   	 		\\
\multicolumn{2}{l}{$ \beta =1 $} \vline &  \multicolumn{2}{c}{$ sign(x) \alpha  |x|^{\alpha-1} $} \\
	  \hline
		\end{tabular}
			\vspace{5pt}
	\caption{Limit behavior of the Fractal variation of $|x|^\alpha, \, \alpha >0$ }
	\label{tab:fracvar}
\end{table}
\begin{figure}[!ht]
	\centering
	 \def\JPicScale{0.8}
	 \ifx\JPicScale\undefined\def\JPicScale{1}\fi
	 \unitlength \JPicScale mm
	 \begin{picture}(80.27,55.96)(0,0)
	 \linethickness{0.3mm}
	 \put(40.21,5.96){\line(0,1){50}}
	 \put(40.21,55.96){\vector(0,1){0.12}}
	 \linethickness{0.3mm}
	 \put(0.21,5.96){\line(1,0){80}}
	 \put(80.21,5.96){\vector(1,0){0.12}}
	 \linethickness{0.5mm}
	 \multiput(40.21,5.96)(0.12,0.12){167}{\line(1,0){0.12}}
	 \linethickness{0.3mm}
	 \put(40.21,25.96){\line(1,0){20.06}}
	 \put(60.21,3.96){\makebox(0,0)[cc]{1.0}}
	 
	 \put(20.21,3.96){\makebox(0,0)[cc]{-1.0}}
	 
	 \put(36.21,29.96){\makebox(0,0)[cc]{1.0}}
	 
	 \linethickness{0.5mm}
	 \multiput(60.21,25.96)(1.38,1.38){15}{\multiput(0,0)(0.11,0.11){6}{\line(0,1){0.11}}}
	 \put(74.68,2.45){\makebox(0,0)[cc]{$\alpha$}}
	 
	 \put(36.21,53.96){\makebox(0,0)[cc]{$\beta$}}
	 
	 \put(28.24,37.93){\makebox(0,0)[cc]{$-\infty$}}
	 
	 \put(28.21,13.96){\makebox(0,0)[cc]{$-\infty$}}
	 
	 \put(46.56,9.83){\makebox(0,0)[cc]{1}}
	 
	 \put(50.27,27.93){\makebox(0,0)[cc]{$g^\prime(x) = \infty$}}
	 
	 \put(78,34){\makebox(0,0)[cc]{0}}
	 
	 \put(46,19){\makebox(0,0)[cc]{$\infty$}}
	 
	 \put(60.21,15.96){\makebox(0,0)[cc]{0}}
	 
	 \put(52.18,37.93){\makebox(0,0)[cc]{$\infty$}}
	 
	 \linethickness{0.3mm}
	 \multiput(0.05,25.96)(1.95,0){21}{\line(1,0){0.98}}
	 \put(78,42){\makebox(0,0)[cc]{0}}
	 
	 \linethickness{0.3mm}
	 \put(60,6){\line(0,1){1}}
	 \linethickness{0.3mm}
	 \put(20,6){\line(0,1){1}}
	 \put(20,28){\makebox(0,0)[cc]{$g^\prime(x) = -\infty$}}
	 
	 \linethickness{0.3mm}
	 \multiput(60.21,25.88)(1.91,0.01){11}{\multiput(0,0)(0.95,0){1}{\line(1,0){0.95}}}
	 \end{picture}
	\caption{Limit behavior, $\upsilon^{\epsilon+}_{\beta} [x^\alpha ], \ x=0$}
	\label{fig:scaling}
\end{figure}

Let $g(x)= |x|^{-\alpha}, \ \alpha >0$. 
Then the limiting behavior of the Fractal variation $\llim {\epsilon} {0} {\upsilon_{\beta}^{ \epsilon  } [g] (x)} $ can  be summarized in   Table \ref{tab:fracvar1}:
\begin{table}[!ht]
	\centering
		\begin{tabular}[t]{l|r|r}
		\hline
				-   	& $x=0$   	& 	 	$ |x| >0 $  \\
	 \hline
	$ \beta <1$     	& $-\infty$					&			 0     	 		\\
		\hline
	$ \beta >1 $ &  \multicolumn{2}{c}{$ -\infty$}\\
  $ \beta =1 $ &  \multicolumn{2}{c}{$  -sign(x) \alpha  |x|^{-\alpha-1} $} \\
	  \hline
		\end{tabular}
			\vspace{5pt}
	\caption{Limit behavior of the Fractal variation of $|x|^{-\alpha}, \, \alpha >0$ }
	\label{tab:fracvar1}
\end{table}
\begin{figure}[!ht]
	\centering
		 \def\JPicScale{0.8}
		 \ifx\JPicScale\undefined\def\JPicScale{1}\fi
		 \unitlength \JPicScale mm
		 \begin{picture}(80.21,55.96)(0,0)
		 \linethickness{0.3mm}
		 \put(40.21,5.96){\line(0,1){50}}
		 \put(40.21,55.96){\vector(0,1){0.12}}
		 \linethickness{0.3mm}
		 \put(0.21,5.96){\line(1,0){80}}
		 \put(80.21,5.96){\vector(1,0){0.12}}
		 \linethickness{0.3mm}
		 \put(0.21,25.96){\line(1,0){80}}
		 \put(60.21,3.96){\makebox(0,0)[cc]{1.0}}
		 
		 \put(20.21,3.96){\makebox(0,0)[cc]{-1.0}}
		 
		 \put(36.21,29.96){\makebox(0,0)[cc]{1.0}}
		 
		 \put(74.68,2.45){\makebox(0,0)[cc]{$\alpha$}}
		 
		 \put(36.21,53.96){\makebox(0,0)[cc]{$\beta$}}
		 
		 \put(19.89,38.03){\makebox(0,0)[cc]{$-\infty$}}
		 
		 \put(50.27,27.93){\makebox(0,0)[cc]{$g^\prime(x)$}}
		 
		 \put(20.11,16.12){\makebox(0,0)[cc]{0}}
		 
		 \put(60.21,15.96){\makebox(0,0)[cc]{0}}
		 
		 \put(60.69,37.93){\makebox(0,0)[cc]{$\infty$}}
		 
		 \linethickness{0.3mm}
		 \put(59.77,6.13){\line(0,1){1}}
		 \linethickness{0.3mm}
		 \put(20,6){\line(0,1){1}}
		 \end{picture}
	\caption{Limit behavior, $\upsilon^{\epsilon+}_{\beta} [x^\alpha ], \ |x|>0$}
	\label{fig:scaling1}
\end{figure}
\end{theorem}

The limiting behavior is represented graphically in Figs. \ref{fig:scaling} and \ref{fig:scaling1}.

\subsection{Fractal variation of smooth functions}
\label{sec:smooth}
We will prove a general theorem allowing one to compute the limit of the fractal variation for smooth functions.
\begin{theorem}[Limit of Fractal Variation about a point]
	\label{th:mapping deriv}
	Let $f(x) \in \fclass{C}{1}$. Then
	\begin{align}
	\llim{\epsilon}{0}{\fracvarplus{f}{x}{\beta}} = \frac{1}{\beta}	\llim{\epsilon}{0}{} \epsilon^{1-\beta} f^{\prime} (x + \epsilon) \\
	\llim{\epsilon}{0}{\fracvarmin{f}{x}{\beta}} = \frac{1}{\beta}	\llim{\epsilon}{0}{} \epsilon^{1-\beta} f^{\prime} (x - \epsilon) 
	\end{align}
\end{theorem}

\begin{corollary}[Vanishing variation theorem]
	\label{th:diffvar1}
	Let $f(x) \in \fclass{C}{1}$ about $x$ and $0<\beta <1$ then
	\[
	\llim {\epsilon} {0} { \fracvarplus{f}{x_0}{\beta} } = 0
	\]
\end{corollary}
This result corresponds with the result obtained in \cite{Chen2010}.

\subsection{Fractal Variation of functions with singular derivatives}
\label{sec:singular}
Closely related to the results of the Vanishing Variation Theorem is the next Theorem.
First let's introduce the concept of the \textit{critical exponent} acting around a singularity of a function.
\begin{definition}
Let $g(x)$ be continuous  function having a singularity at $x_s$.

Then the left critical exponent $\alpha$ be the minimal exponent in the power term for which the quantity
$ h^\alpha \, \mathcal{T}_{-h}[g](x) $ is finite at $x_s$.  
Or formally,
\[
\mathcal{P}_{+}[g](x=x_s):=\alpha \ \left| \left\{  0< \inf\limits_{\alpha} \left|{ \llim{h}{0}{h^\alpha \, \mathcal{T}_{-h}[g](x)}}  \right| < \infty \right\} \right.
\]
where $h>0$.
Let the right critical exponent $\alpha$ be the minimal exponent in the power term for which 
\[
\mathcal{P}_{-}[g](x=x_s):=  \alpha \ \left| \left\{ 0<  \inf\limits_{\alpha} \left|{ \llim{h}{0}{h^\alpha \, \mathcal{T}_{h}[g](x)}}  \right| < \infty \right\} \right.
\] 
When $g(x)$ is bounded $\mathcal{P}[g](x)$ will be assumed 0.
\end{definition}
\begin{theorem}[Singular variation theorem]
\label{th:diffvar2}
Let $f(x) \in \fclass{C}{n}$ has a singular derivative at $x_s$ 
$f^\prime(x)$ be such that 
$\mathcal{P}_{+}[f^\prime](x_s)=\alpha$.
Then in limit
\begin{description}
	\item[ for $ \beta < | 1 - \alpha|$] $\llim{\epsilon}{0}{ } \fracvarplus{f}{x_s -}{\beta} =0$
		\item[ for $ \beta = | 1 - \alpha|$] $\llim{\epsilon}{0}{ } \fracvarplus{f}{x_s -}{\beta} $ is finite
			\item[ for $ \beta > | 1 - \alpha|$] $\llim{\epsilon}{0}{ } \fracvarplus{f}{x_s-}{\beta}$ is unbounded.	
\end{description}
	assuming always $\beta \in [0 \ 1)$.
	For $x \neq x_s $ $\llim{\epsilon}{0}{ } \fracvarplus{f}{x_s -}{\beta} =0$. 
	The result can be illustrated from Figs. \ref{fig:scaling} and \ref{fig:scaling1}.
\end{theorem}
\begin{theorem}
\label{th:diffvar3}
Let $f(x) \in \fclass{C}{n}$ is such that 
$\mathcal{P}_{-}[f^{\prime}](x_s)=\alpha$. Then in limit
\begin{description}
	\item[ for $ \beta < | 1 - \alpha|$] $\llim{\epsilon}{0}{ } \fracvarmin{f}{x_s +}{\beta} =0$
		\item[ for $ \beta = | 1 - \alpha|$] $\llim{\epsilon}{0}{ } \fracvarmin{f}{x_s +}{\beta} $ is finite
			\item[ for $ \beta > | 1 - \alpha|$] $\llim{\epsilon}{0}{ } \fracvarmin{f}{x_s+}{\beta} $ is unbounded.	
\end{description}
	assuming always $\beta \in [0 \ 1)$.
	The proof is analogous to the proof of Theorem \ref{th:diffvar2}.
For $x \neq x_s $ $\llim{\epsilon}{0}{ } \fracvarplus{f}{x_s +}{\beta} =0$. 
\end{theorem}

\begin{definition}
Let $u (x) $ be the function with the following definition:
\[
u(x):=\llim {n}{0}{u_n (x)}
\]
defined as the limiting sequence of 
\[
u_n (x):= \left\{
\begin{array}{ll}
0 ,& x <0 \\   
1 ,&  x < \epsilon_n \\
0 ,& x\geq \epsilon_n
\end{array}
\right.
\]
for $\epsilon_n = \frac{\upsilon}{2^n}$ where $0<\upsilon <2$ is an arbitrary small number, $\ n \in \fclass{N}{}$.
\end{definition}
\begin{example}
The result can easily be extended for  scaled and translated versions of $x^\alpha$ :
$f(x) = (s \,x-x_0)^\alpha$ because of the commutativity of translation and the homogeneity property. Then 
\begin{equation}
\llim{\epsilon}{0}\upsilon^{\epsilon+}_\alpha \left[ \, |x -x_0|^\alpha \right] = \mathbf{u}(x - x_0)
\end{equation}
for $1>\alpha>0$.
\end{example}

\subsection{Fractal variation of functions with singularities}
\label{sec:fvarsingular}

The behavior of the fractal variation is especially interesting when the argument function can become singular. It can be demonstrated that the fractal variation preserves singularities in limit. 
In the standard setting, the next results have to be interpreted in distributional sense. 

 \subsection{Fractal variation of $\delta$-sequences}
 \label{sec:delta}

 The difficulty in the standard treatment of the Delta function's properties comes from the fact that the function has to assume unlimited growth at the origin and therefore is discontinuous. 
 This makes the notion difficult to handle in the standard framework of analysis \cite{Todorov1990} (recent review on historical developments in \cite{Katz2013}). 
 With standard analytic arguments the notion of "Delta" can be handled in Swartz distribution theory.
 In the following we give elementary treatment of some of the  properties of the Dirac's Delta function. 
 We will always assume that the value at the origin can be defined in some sense. 

\begin{figure}[ht!]
	\centering
 
		\def\JPicScale{1.2}
	 
	 \ifx\JPicScale\undefined\def\JPicScale{1}\fi
	 \unitlength \JPicScale mm
	 \begin{picture}(54.89,40)(0,0)
	 \linethickness{0.3mm}
	 \put(4.89,4.79){\line(1,0){50}}
	 \put(54.89,4.79){\vector(1,0){0.12}}
	 \linethickness{0.3mm}
	 \put(30,5){\line(0,1){35}}
	 \put(30,40){\vector(0,1){0.12}}
	 \linethickness{0.3mm}
	 \put(25.1,3.74){\line(0,1){1}}
	 \put(39.91,27.6){\makebox(0,0)[cr]{$\frac{1}{ \epsilon_n }$}}
	 
	 \put(25.32,-0.16){\makebox(0,0)[cc]{$-\epsilon_n \over 2$}}
	 
	 \put(35.27,-0.27){\makebox(0,0)[cc]{$\epsilon_n \over 2$}}
	 
	 \linethickness{0.3mm}
	 \put(20.24,14.95){\line(1,0){19.73}}
	 \multiput(20.21,4.8)(0.03,10.15){1}{\line(0,1){10.15}}
	 \multiput(39.94,4.8)(0.03,10.15){1}{\line(0,1){10.15}}
	 \put(20.21,4.8){\line(1,0){19.73}}
	 \put(44.95,17.35){\makebox(0,0)[cr]{$\frac{1}{2 \epsilon_n }$}}
	 
	 \linethickness{0.3mm}
	 \multiput(25.06,24.95)(2.17,0){5}{\line(1,0){1.09}}
	 \multiput(25.06,24.95)(0,-1.92){11}{\multiput(0,0)(0,-0.96){1}{\line(0,-1){0.96}}}
	 \multiput(34.84,24.95)(0,-1.92){11}{\multiput(0,0)(0,-0.96){1}{\line(0,-1){0.96}}}
	 \multiput(25.09,4.79)(2.17,0){5}{\line(1,0){1.09}}
	 \linethickness{0.3mm}
	 \put(20.21,3.75){\line(0,1){1}}
	 \put(20.27,2.23){\makebox(0,0)[cc]{$-\epsilon_n$}}
	 
	 \put(40.04,2.24){\makebox(0,0)[cc]{$ \epsilon_n$}}
	 
	 \linethickness{0.3mm}
	 \put(34.88,3.8){\line(0,1){1}}
	 \linethickness{0.3mm}
	 \put(39.95,3.81){\line(0,1){1}}
	 \linethickness{0.3mm}
	 \multiput(25.11,4.73)(0.47,1.93){11}{\multiput(0,0)(0.12,0.48){2}{\line(0,1){0.48}}}
	 \linethickness{0.3mm}
	 \multiput(30.05,24.89)(0.46,-1.91){11}{\multiput(0,0)(0.11,-0.48){2}{\line(0,-1){0.48}}}
	 \linethickness{0.3mm}
	 \multiput(20.21,4.84)(0.12,0.12){81}{\line(0,1){0.12}}
	 \linethickness{0.3mm}
	 \multiput(29.89,14.95)(0.12,-0.12){84}{\line(0,-1){0.12}}
	 \put(2.18,27.07){\makebox(0,0)[cc]{}}
	 
	 \put(15.32,27.07){\makebox(0,0)[cc]{}}
	 
	 \end{picture}
 
	\caption{Graph of two Delta sequences}
	\label{fig:DeltaSequence}
The graph of the triangular sequence is rescaled by factor 2.
\end{figure}
\begin{theorem}[Fractal variation of the $\delta$-function]
	\label{th:vardelta}
Let $\delta(x)$ be the Dirac's delta function/distribution 
\[
\delta(x):=\llim {n}{\infty}{\delta_n (x)}
\]
defined as the limiting sequence of functions-prototypes (i.e. predistributions)
 \[
 \delta_n(x):= \frac{1}{s_n} \ \psi\left(\frac{x}{s_n} \right)
 \]
 having the following  properties:
 \begin{itemize}
 	\item $ \int\limits_{-\infty}^{\infty} \psi\left( x \right) dx =1 $
 	\item $\psi\left( x   \right) = \psi\left( -x  \right) $
 	\item positive at the origin 
 	\item monotonously decreasing from the origin towards $\pm \infty$  as $\sim \frac{1}{x^2}$.
 \end{itemize}
 that are 
  parametrized by the Cauchy sequence $\llim{n}{\infty}{s_n =0}$ with $s_1 \leq 1$ 
Then formally in limit
\begin{equation}
\label{eq:deltafrac}
\llim{\epsilon}{0} {\fracvarplus {\delta} {x}{\beta}} = -\frac{sign(x)}{|x|^\beta}\; \delta(x)
\end{equation}
\end{theorem}
We are going to prove two lemmas serving as limiting cases.
\begin{lemma}[Fractal variation of the pulsed $\delta$-sequence]
	\label{th:vardeltastep}
Let $\delta(x)$ be a rectangular sequence
 \[
\delta_n (x):= \left\{
	\begin{array}{ll}   
		0  ,& |x| \geq \epsilon_n \\      
	\frac{1}{2 \epsilon_n }  ,&  |x| < \epsilon_n 
	\end{array}
\right.
\]
where $\epsilon_n = \frac{\upsilon}{2^n}$, $0<\upsilon <2$ is an arbitrary small number, $\ n \in \fclass{N}{}$.

Then for $\beta \leq 1$
\[
   \fracvar {\delta}{x}{\beta}{\epsilon_n+}   = \left\{
	\begin{array}{ll}   
		0  ,& |x| \geq \epsilon_n \\      
	- \frac{ sign(x)}{2 \epsilon_n^{\beta +1} }  ,&  \epsilon_{n-1} \leq |x| < \epsilon_n \\
		0  ,& |x| < \epsilon_{n-1} \\   
	\end{array}
\right.
\]

\end{lemma}
%
%
\begin{lemma}[Fractal variation of the triangular $\delta$-sequence]
	\label{th:vardelta1piece}
	Let  $\delta_n (x)$ be defined as a sequence of triangular functions
	\[
	\delta_n (x):= \left\{
	\begin{array}{ll}   
	0  ,& |x| \geq \epsilon_n \\      
	\frac{1}{\epsilon_n } - \frac{ x}{\epsilon_n^2 } ,&  0 \leq x < \epsilon_n \\
	\frac{1}{\epsilon_n } + \frac{ x}{\epsilon_n^2 } ,&  0 \geq x > -\epsilon_n

	\end{array}
	\right.
	\]
	where $	\epsilon_n = \frac{\upsilon}{2^n}$, $0<\upsilon <2$ is an arbitrary small number, $\ n \in \fclass{N}{}$.
	
	Then for $\beta \leq 1$
	\[
	\fracvar {\delta_n}{x}{\beta}{\epsilon_n+}   = \left\{
	\begin{array}{ll}   
	0  ,& |x| \geq \epsilon_n \\      
	- \frac{ sign(x)}{2 \epsilon_n^{\beta +1} }  ,&  \epsilon_{n-1} \leq |x| < \epsilon_n \\
	0  ,& |x| < \epsilon_{n-1} \\   
	\end{array}
	\right.
	\]

\end{lemma}
%
Therefore,  any scale-dependent parametric smooth function constrained between the rectangular pulse and the triangular pulse at some scales will also exhibit this scaling limit behavior behavior.

  Central argument for the subsequent presentation will be differentiating property of the $\delta$-function:
  \[
  \int_{-\infty}^{+\infty} \delta'(x)f(x) \,dx =\left.\delta(x)f(x)\;\right|_{-\infty}^{+\infty} -\int_{-\infty}^{+\infty}  \delta(x)f'(x)\,dx =-f'(0)
  \]
  Indeed by arguments of integration by parts for the prototype we have
    \[
    \int_{-\infty}^{+\infty}  \psi'(x)f(x) \,dx =\left.\psi(x)f(x)\right|_{-\infty}^{+\infty} -\int_{-\infty}^{+\infty}  \psi(x)f'(x)\,dx  
    \]
  Then by mapping to scale
    \[
    \int_{-\infty}^{+\infty}  \frac{1}{s_n^2} \, \psi' \left({x \over s_n} \right) f(x) \,dx =
    \frac{1}{s_n^2} \left.\psi\left( {x \over s_n} \right)f(x)\right|_{-\infty}^{+\infty}  - \int_{-\infty}^\infty \frac{1}{s_n^2} \, \psi \left({x \over s_n} \right) f'(x)\,dx  
    \]

  Let $f_n(x) \in \fclass{C}{\infty}$ be a function of a Delta sequence. That is
  \[
  f_n(x)= \frac{1}{s_n} \ \psi\left(\frac{x}{s_n} \right)
  \]
 Due to the symmetry about the origin it can be demonstrated that $\psi^\prime\left( x  \right) = - \psi^\prime\left( -x  \right) $,  $\psi^{\prime \prime}\left( x  \right) =  \psi^{\prime\prime}\left( -x  \right) $ and so on.
 From these properties it follows in particular that $ \psi^\prime\left( 0 \right) =0 $ and $ \psi^{\mathrm{iii}}\left( 0 \right) =0 $.

 $f(x)$ can be expanded in a  4\textsuperscript{th} order Taylor series about the origin as
 \[
 f(x) =\frac{1}{s} \left( \psi\left( 0  \right)  + \frac{1}{2\, s^2} \psi^{\prime \prime}(0) \, x^2  +
 \frac{1}{24 \, s^4}  \psi^{\mathrm{i v}}(0) \, x^4 + \mathcal{O}(x^6)
 \right)
 \]
 and
 \[
 f^\prime(x) = \frac{1}{s} \left(  \frac{1}{s^2} \psi^{\prime \prime}(0) \, x  +
 \frac{1}{6\, s^4}  \psi^{\mathrm{i v}}(0) \, x^3 + \mathcal{O}(x^5)  \right)
 \]
 From the area condition it follows that $\llim{x}{\pm\infty}{\psi(x) =0}$.
 Since for large $x$ $\psi(x)$ must decay very fast towards \textit{0} the signs of its derivatives must alternate in order to have cancellation. 
 Let's make the substitution $ \psi^{\prime \prime}(0) = b$ and  $\psi^{\mathrm{i v}}(0)= c$. 
 Then by assumption $ b<0$.
 Then if follows that about the origin
 \[
 f^{\prime}(x) = \frac{b\,x}{{s}^{3}} + \frac{c\,{x}^{3}}{6\,{s}^{5}}
 \]

 which has  extrema at
 \[
 x_m = \pm  \sqrt{-\frac{2 \ b}{ c}} \, s
 \]
 Since  \textit{b} and \textit{c} have different signs the roots $x_m$ are all real. 
 \[
 f^\prime(x_m) =
 \frac{ 2\, b}{3\,{s}^{2}} \,\sqrt{-\frac{2\, b}{c}}
 \]
 
 On the other hand, if we consider the  4\textsuperscript{nd} order Taylor expansion, the fractal variation at the positive extremum is
 \[
  \fracvarplus{f}{x_m}{\beta} = - \frac{{\epsilon}^{1- \beta}}{s^2 \, 24\,{\sqrt{2}}}\,  \left(  16\,  \left( \sqrt{-\frac{b}{c}}  \right)^{3}\,c + 48\,b\,\sqrt{-\frac{b}{c}}  \, +8\,\sqrt{-\frac{b}{c}}\,c\, \left( \frac{\epsilon}{s}\right) ^{2}+\sqrt{2}\,c\,
  	\left( \frac{\epsilon }{s}\right) ^{3}\right)  
\]
 Theretofore, the behavior of the extremum is dictated by two exponents:
 $ \frac{ {\epsilon}^{1-\beta} }{{s}^{ 2}}$ and $\frac{\epsilon}{s}$.
 Let's make the substitution
 $ s= k \, \epsilon^p $. 
 Then
 \[
 \fracvarplus{f}{x_m}{\beta} =  
 -\frac{2\sqrt{2}\,b\,\sqrt{-\frac{b}{c}}\,{\epsilon}^{-2\,p-\beta+1}}{3\,{k}^{2}}-\frac{\sqrt{-\frac{b}{c}}\,c\,{\epsilon}^{-4\,p-\beta+3}}{3\,\sqrt{2}\,{k}^{4}}-\frac{c\,{\epsilon}^{-5\,p-\beta+4}}{24\,{k}^{5}}
 \]
 Then if we wish to retain the singularity while decreasing scale:
\[
   p > \frac{1-\beta}{2}  \cap    p > \frac{3-\beta}{4} \cap  p > \frac{4-\beta}{5}  
\]
 Therefore, only $ p > \frac{1-\beta}{2} $
is an admissible exponent in the substitution $ s=  \epsilon^p $. 
 
 This amounts to a scale-relative $s/\epsilon$-limit procedure. 
 
 In particular, if $p=1$ and $k=1$ then
 \[
 \fracvarplus{f}{x_m}{\beta} =  
 -\frac{\left(  16\,  \left( \sqrt{-\frac{b}{c}}  \right)^{3}+8\,\sqrt{-\frac{b}{c}}+\sqrt{2}\right) \,c
 	+48\,b\,\sqrt{-\frac{b}{c}}}{24\,{\sqrt{2}}\, {\epsilon}^{1+ \beta}}
 \]
 which scales as prescribed in Th. \ref{th:vardelta}.
 
 There is a critical value of the order $\beta=\frac{1}{3}$ for which separation of limits can not be used and power-law substitution can not be applied.
 The result is demonstrated in Fig. \ref{fig:lim}.
 \begin{figure}[!ht]
 	\centering
 	\ifx\JPicScale\undefined\def\JPicScale{1}\fi
 	\unitlength \JPicScale mm
 	\begin{picture}(60.69,59.45)(0,0)
 	\linethickness{0.3mm}
 	\multiput(5,27.45)(0.24,-0.12){208}{\line(1,0){0.24}}
 	\linethickness{0.3mm}
 	\multiput(5,2.45)(0.12,0.12){417}{\line(1,0){0.12}}
 	\linethickness{0.3mm}
 	\put(5.05,2.45){\line(1,0){55.64}}
 	\put(60.69,2.45){\vector(1,0){0.12}}
 	\linethickness{0.3mm}
 	\put(5,2.45){\line(0,1){57}}
 	\put(5,59.45){\vector(0,1){0.12}}
 	\linethickness{0.3mm}
 	\put(21.68,1.45){\line(0,1){1}}
 	\linethickness{0.3mm}
 	\put(30,1.45){\line(0,1){1}}
 	\linethickness{0.3mm}
 	\put(4,27.45){\line(1,0){1}}
 	\linethickness{0.3mm}
 	\multiput(5,52.45)(1.96,0){26}{\line(1,0){0.98}}
 	\linethickness{0.3mm}
 	\multiput(55,2.45)(0,1.96){26}{\line(0,1){0.98}}
 	\linethickness{0.3mm}
 	\multiput(21.65,2.5)(0,1.96){9}{\line(0,1){0.98}}
 	\linethickness{0.3mm}
 	\put(4,52.45){\line(1,0){1}}
 	\put(2,27.45){\makebox(0,0)[cc]{$\frac{1}{2}$}}
 	
 	\put(2,52.45){\makebox(0,0)[cc]{1}}
 	
 	\put(60.21,-0.9){\makebox(0,0)[cc]{$\beta$}}
 	
 	\put(2.21,57.86){\makebox(0,0)[cc]{p}}
 	
 	\put(21.7,-1.12){\makebox(0,0)[cc]{$\frac{1}{3}$}}
 	
 	\put(29.95,9.47){\makebox(0,0)[cc]{$0$}}
 	
 	\put(30.05,34.15){\makebox(0,0)[cc]{$\delta(x)$}}
 	
 	\put(36.49,24.95){\makebox(0,0)[cc]{}}
 	
 	\put(41.54,11.86){\makebox(0,0)[cc]{$u(x)$}}
 	
 	\put(54.89,-0.96){\makebox(0,0)[cc]{1}}
 	
 	\put(4.79,-0.8){\makebox(0,0)[cc]{0}}
 	
 	\put(2.23,3.24){\makebox(0,0)[cc]{0}}
 	
 	\put(30.04,-1.05){\makebox(0,0)[cc]{$\frac{1}{2}$}}
 	
 	\end{picture}
 	\caption{Critical exponents of the $s/\epsilon$-limit procedure}
 	The primitive function maps to $u(x)$ along the critical line $1- 2p - b=0$.
 	\label{fig:lim}
 \end{figure}
 Only in the limit $\beta=1$, which corresponds to usual differentiation the $\epsilon$ and $s$-limiting procedures seem to be unrelated.
 
 This principle corresponds with the theory of scale relativity where the scale $\epsilon$-resolution  is related to the temporal $t$-resolution as argued heuristically by Nottale \cite{Nottale2013}. 
 In other words, not all trajectories towards singularity in the resolution space are admissible if the separation of limits is used.

 \begin{corollary}
 	\label{corr:deltadfracvar}
 	Formally we can write the following equation for the limit of the fractal variation in distributional sense:
 	\begin{equation}
 	\llim{\epsilon}{0}{\fracvarplus{f}{x}{\beta}}=  -\int\limits_{-\infty}^{+\infty}\frac{\delta (x-z)}{|x-z|^\beta} f(z)\, dz
 	\end{equation}
 	where $\beta \leq 1$.
 \end{corollary}
 
 \begin{theorem}[Singular variation theorem ii]
 	\label{th:diffvar3}
 	Let $f(x) \in \fclass{C}{n}$ have a singularity at $x_s$ such that 	$\mathcal{P}_{+}[f ](x_s)=\alpha$ and 	$\mathcal{P}_{-}[f ](x_s)=\alpha$.
 	Then in limit
 	\[
 	\llim{\epsilon} {0} { \fracvarplus{f}{x}{\beta}} = K \delta_{ \alpha + \beta -1 } \left( x- x_s\right) 
 	\]
 	for some real \textit{K}.
 \end{theorem}
 \begin{example}
 	\label{prop:fractan}
 	For    $ \beta < 1$,
 	\[
 	\llim{\epsilon} {0} {\upsilon_{\beta}^{ \epsilon +} [\tan]}(x)  = \sum\limits^{\infty}_{k = 0}{ \delta \left(x- \frac{(2 k +1) \pi}{2} \right)}
 	\] where $k \in \fclass{Z}{}$.
 	
 	According to Theorem \ref{th:diffvar1} when $x \neq \frac{\pi}{2}$ then 
 	$ \llim{\epsilon} {0} {\upsilon_{\beta}^{ \epsilon +} [\tan]}(x)  = 0
 	$. When $x \rightarrow \frac{\pi}{2}$
 	\[\llim{x, \epsilon} {0} {  
 		\frac{\tan \left(x + \frac{\pi}{2} + \epsilon \right) - \tan \left(x + \frac{\pi}{2} \right)}{\epsilon^\beta}} =
 	\llim{x, \epsilon} {0} {}
 	-\frac{2\,\sin\left( \epsilon\right) }{\left(\cos {(2\,x+\epsilon)} -\cos{( \epsilon)} \right) \epsilon^\beta}
 	\]
 	We make the anasatz $x=\epsilon$
 	\[
 	\llim{ \epsilon} {0} {-\frac{2\,\mathrm{sin}\left( \epsilon\right) }{\left( \cos{(3\,\epsilon)} -\cos{( \epsilon)}\right) \epsilon^\beta}} =
 	\llim{ \epsilon} {0} {
 		\frac{1}{\sin\left( 2\,\epsilon\right)\epsilon^\beta }
 	} = +\infty
 	\]
 \end{example}
 Similar arguments lead to the next result
 \begin{example}
 	\label{prop:fraccot}
 	For    $ \beta < 1$,
 	\[
 	\llim{\epsilon} {0} {\upsilon_{\beta}^{ \epsilon +} [\cot]}(x)  = \sum\limits^{\infty}_{k = 0}{ \delta \left(x- ( k +1) \pi \right)}
 	\] where $k \in \fclass{Z}{}$.
 \end{example}
 

  
 \begin{definition}
 	\label{def:deltamap}
  Let us define the parametric scaling  map $S_a:   \fclass{C}{0} \rightarrow \fclass{C}{0} $ for functions $f: \fclass{R}{} \rightarrow \fclass{R}{} $ with the following properties
	\begin{align}
 S_a:  \  f(x)  &\rightarrow \, f(a \, x) \\
 S_a:  \  \int\limits_{x}^{y} f(t) \,dt &\rightarrow \int\limits_{x \, a}^{y \, a} f(t) \,dt
	\end{align}
where we interpret the symbols $x$ and $y$ as external variables and $a>0$.
The $n$-fold composition of operations will be denoted by $S_{a}^{n}$.
Next, we define also the \textit{delta} map $\hat{\delta}_a:= a \circ S_a$.
\end{definition}
\begin{proposition}
	From the definitions is follows that
	\[
	S_a \circ \int = a  \int \circ \ S_a = \int \circ \ \hat{\delta}_a	
	\]
\end{proposition}
\begin{proposition}
	\label{prop:scaling1}
Let $a>1$ and $f(x)$ is such that $\llim{x}{\infty}{f(x) = 0}$ then 
\[
\llim{n}{\infty}{(S_a)^n} f(x) = \left\{
\begin{array}{ll}   
f(0)  ,& x=0 \\    
0  ,& |x| \geq 0 \\      
\end{array}
\right.
\]
\end{proposition}
 \begin{proposition}
 	Let $\psi (x) \in \fclass{C}{\infty}$ has the following additional properties:
 	\begin{itemize}
 		\item $ \int\limits_{-\infty}^{\infty} \psi\left( x \right) dx =1 $
 		\item monotonously decreasing in both directions from the origin.
 	\end{itemize}
 	Then we have that the following 
 	$\llim{n}{\infty}   (\hat{\delta}_a)^n \;\psi (x)$ exists.
 	Moreover,
	\[
	{\int\limits_{x}^{y}}  \llim{n}{\infty}{} \hat{\delta}_a^n \;\psi(x) \; dx = \left\{
	\begin{array}{ll}   
	0  , & 0 \notin [x, y] \\      
	1, &  0 \in [x, y] \\
	\end{array}
	\right.
	\]
 if $a>1$.
 \end{proposition}
%
\begin{proposition}
Let the scale derivative be defined as the limit of the symmetrical central difference:
\[
\hat{\psi}^\prime(x):=\llim{a}{0}{} \frac{  \psi \left({{z+a^2/2}\over{a}}\right) 
	-  \psi \left({{z-a^2/2}\over{a}}\right)  }{a} 
\]
Then,

\[
\hat{\psi}^\prime(x)=\left\{
\begin{array}{ll}   
0  , & x \neq 0\\      
\frac{\psi_{+}^\prime (0) +\psi_{-}^\prime (0)}{2}, &  x =0 \\
\end{array}
\right. 
\]
\end{proposition}

\section*{Acknowledgments}
The work has been supported in part by a grant from Research Fund - Flanders (FWO), contract number 0880.212.840.
The author would like to acknowledge Dr. Fay\c{c}al Ben Adda for helpful feedback.

\bibliographystyle{plain}  
\bibliography{fracvar}

\begin{thebibliography}{10}

\bibitem{Adda1997}
F.~Ben Adda.
\newblock Geometric interpretation of the fractional derivative.
\newblock {\em J. Fract Calc}, 11:21 -- 51, 1997.

\bibitem{Adda2001}
F.~Ben Adda and J.~Cresson.
\newblock About non-differentiable functions.
\newblock {\em J. Math. Analysis Appl.}, 263:721--737, 2001.

\bibitem{Adda2005}
F.~Ben Adda and J.~Cresson.
\newblock Fractional differential equations and the {Schr\"odinger} equation.
\newblock {\em App. Math. Comp.}, 161:324--345, 2005.

\bibitem{Adda2013}
F.~Ben Adda and J.~Cresson.
\newblock Corrigendum to {"About non-differentiable functions" [J. Math. Anal.
  Appl. 263 (2001) 721 -- 737]}.
\newblock {\em J. Math. Analysis Appl.}, 408(1):409 -- 413, 2013.

\bibitem{Babakhani2002}
A.~Babakhani and V.~Daftardar-Gejji.
\newblock On calculus of local fractional derivatives.
\newblock {\em J. Math. Analysis Appl.}, 270(1):66 -- 79, 2002.

\bibitem{Bagley1983}
R.~Bagley and P.~J. Torvik.
\newblock A theoretical basis for the application of fractional calculus to
  viscoelasticity.
\newblock {\em J. Rheology}, 27:201 -- 210, 1983.

\bibitem{Caputo1967}
M.~Caputo.
\newblock Linear models of dissipation whose {Q} is almost frequency
  independent {II}.
\newblock {\em Geophys J R Ast Soc}, 13(529):529 --539, 1967.

\bibitem{Caputo1971}
M.~Caputo and F.~Mainardi.
\newblock Linear models of dissipation in anelastic solids.
\newblock {\em Rivista del Nuovo Cimento}, 1:161 -- 198, 1971.

\bibitem{Chen2010}
Y.~Chen, Y.~Yan, and K.~Zhang.
\newblock On the local fractional derivative.
\newblock {\em J. Math. Anal. Appl.}, pages 17 -- 33, 2010.

\bibitem{Jumarie2006}
G.~Jumarie.
\newblock Modified {Riemann-Liouville} derivative and fractional {Taylor}
  series of non-differentiable functions.
\newblock {\em Comp. Math. Appl.}, 51(9 -- 10):1367 -- 1376, 2006.

\bibitem{Katz2013}
M.G. Katz and D.~Tall.
\newblock A {Cauchy-Dirac} {Delta} function.
\newblock {\em Foundations of Science}, 18(1):107 -- 123, 2013.

\bibitem{Kolwankar2013}
K.~Kolwankar.
\newblock Local fractional calculus: a review.
\newblock Number arXiv:1307.0739v1. 2013.

\bibitem{Kolwankar1997a}
K.~Kolwankar and A.D. Gangal.
\newblock H\"older exponents of irregular signals and local fractional
  derivatives.
\newblock {\em Pramana J. Phys}, 1(1):49 -- 68, 1997.

\bibitem{Mallat1992}
S.~Mallat and W.-L. Hwang.
\newblock Singularity detection and processing with wavelets.
\newblock {\em Information Theory, IEEE Transactions on}, 38(2):617 -- 643,
  1992.

\bibitem{Nottale2013}
L.~Nottale and M.N. C\'el\'erier.
\newblock Emergence of complex and spinor wave functions in scale relativity.
  i. nature of scale variables.
\newblock {\em Journal of Mathematical Physics}, 54(11):--, 2013.

\bibitem{Oldham1970}
K.~Oldham and J.~Spanier.
\newblock {\em The Fractional Calculus}.
\newblock Academic Press, London, 1970.

\bibitem{Ross1977}
B.~Ross.
\newblock The development of fractional calculus 1695 -- 1900.
\newblock {\em Historia Math.}, 4:75 --89, 1977.

\bibitem{Todorov1990}
T.~Todorov.
\newblock A non-standard {Delta} function.
\newblock {\em Proc. Am. Math. Soc}, 110:1143 -- 1144, 1990.

\end{thebibliography}

\end{document}